\documentclass[journal]{IEEEtran}
\usepackage{mathrsfs}
\usepackage{graphicx}
\usepackage{epstopdf}

\usepackage{subeqnarray}
\usepackage{cases}

\usepackage{slashbox}
\usepackage{amssymb} 
\usepackage{amsthm}
\usepackage[numbers,sort&compress]{natbib}
\newtheorem{theorem}{Theorem}
\newtheorem{remark}{Remark}
\newtheorem{example}{Example}
\newtheorem{lemma}{Lemma}
\newtheorem{definition}{Definition}

\ifCLASSINFOpdf

\else
 \fi
\begin{document}
%
\title{Finite-time and Fixed-time Convergence in Continuous-time Optimization}

\author{Yuquan~Chen,~
        Yiheng~Wei,~
		and~YangQuan~Chen,~\IEEEmembership{Senior Member,~IEEE,}
\thanks{Yuquan Chen is with Department of Automation, Hohai University, Nanjing, 210024, P.R. China.
        {\tt\small cyq@mail.ustc.edu.cn}}
\thanks{Yiheng Wei is with School of Mathematics, Southeast University, Nanjing 210096, P.R. China.
        {\tt\small neudawei@seu.edu.cn}}
\thanks{YangQuan Chen is with School of Engineering, University of California, Merced, 5200 North Lake Road, Merced, CA 95343, USA
        {\tt\small ychen53@ucmerced.edu}}
}


\maketitle

\begin{abstract}
It is known that the gradient method can be viewed as a dynamic system where various iterative schemes can be designed as a part of the closed loop system with desirable properties. In this paper, the finite-time and fixed-time convergence in continuous-time optimization are mainly considered. By the advantage of sliding mode control, a finite-time gradient method is proposed, whose convergence time is dependent on initial conditions. To make the convergence time robust to initial conditions, two different designs of fixed-time gradient methods are then provided. One is designed using the property of sine function, whose convergence time is dependent on the frequency of a sine function. The other one is designed using the property of Mittag-Leffler function, whose convergence time is determined by the first positive zero of a Mittag-Leffler function. All the results are extended to more general cases and finally demonstrated by some dedicated simulation examples.
\end{abstract}

\IEEEpeerreviewmaketitle

\section{Introduction}
\IEEEPARstart{G}{radient} method (GM), a classical optimization algorithm, has been widely used in many engineering applications like adaptive filter \citep{song2019affine,qian2020maximum}, artificial intelligence \cite{2018Scalable,2020Gradient}, and system identification \citep{angulo2015nonlinear,lin2015parameter,2017System}. To improve the convergence speed of GM, many variants have been proposed. On one hand, second-order GM, namely Newton-method, is proposed, where Hessian matrix is used to modify the update direction for a faster convergence speed \cite{akgul2019fractional}. However, Hessian matrix is always difficult for computing and quasi-Newton method is then proposed. On the other hand, scholars aim at designing efficient iterative policy to improve the convergence speed, such as GM with momentum \cite{qian1999momentum} and Neserov's accelerated GM (NAGM) \cite{su2014differential}. Moreover, conventional GM is a local algorithm, which can easily get trapped into a local minimum or a saddle point. Perturbation is introduced to help escaping the saddle points more efficiently \cite{jin2017escape} and L{\'e}vy perturbation is then proven to further improve the global convergence capability of GM due to the frequently large jumps \cite{simsekli2019tail}.

Generally, optimization algorithm and system control are closely related to each other. Many iterative algorithms can be formulated as a closed-loop system and the properties can be derived using system analyses methods. \cite{bhaya2003iterative} derives the Newton-Raphson algorithm and conjugate GM by selecting suitable Lyapunov functions. \cite{kashima2007system} views the numerical algorithm as a nonlinear feedback system, and the convergence property is then proven by positive real theorem. Furthermore, by the advantage of control theory, \cite{noroozi2009finite} designs a novel quasi-Newton algorithm, resulting in a finite-time convergence. NAGM has played an important role in machine learning. In \cite{su2014differential}, NAGM is formulated as a second-order differential equation, and its properties are analyzed by Lyapunov theorem. On this basis, more results have been published for analyzing and designing accelerated GMs \cite{wibisono2016variational,wilson2016lyapunov,Laborde2019lyapunov}. Recently, with the fast development of fractional calculus, fractional order design has gained increasing concentrations. By replacing the first-order difference operator with fractional-order difference operator, fractional iterative LMS algorithm is proposed in \cite{tan2015novel}, and it is stated that a larger iterative order resulted in a faster convergence speed while a smaller iterative order resulted in a smaller steady error. A novel fractional LMS algorithm with hybrid iterative order is then proposed in \cite{cheng2017universal} to get both a faster convergence speed and a smaller steady error. A general design for fractional GMs is concluded in \cite{chen2019unified}.

Finite-time convergence is always expected for optimization algorithms rather than asymptotic convergence, which has been discussed in system control for a long time. To achieve finite-time convergence in system control, many efficient methods have been proposed. Among all the finite-time control methods, sliding mode control is the most popular one, where the sign function is used to guarantee the finite-time convergence to the sliding manifold \cite{mishra2018arbitrary}. Though many invariants have been proposed to improve the convergence speed, the convergence time is usually dependent on initial conditions, and the convergence time can be arbitrarily large, which is undesirable. Therefore, scholars pay more attention on fixed-time convergence, where there is an upper bound for the convergence time with arbitrary initial conditions. In\cite{mishra2018arbitrary}, a fixed-time sliding mode controller is designed to realize a fixed-time convergence to the sliding manifold. A basic Lyapunov theorem for fixed-time convergence is then proposed in \cite{polyakov2011nonlinear}. On this basis, many fixed-time controllers are then designed \cite{tian2018fixed,basin2019finite,zuo2018overview}. In \cite{chen20182}, two novel fixed-time reaching laws are proposed from a different perspective, using the property of sine function and Mittag-Leffler function. Recently, finite-time GM has gained increasing attentions and is considered in many practical engineering \cite{yu2017dynamical,lin2016distributed}. However, the convergence time of the mentioned finite-time GMs are all relevant to initial conditions. In \cite{romero2020finite}, the concept of fixed-time GM is somewhat mentioned, but it requires the exact information of initial conditions, which is always unknown to designers.

Motivated by aforementioned reasons, finite-time and fixed-time GMs are considered in this paper. Similar to the design of sliding mode controller, a finite-time GM is designed using the concept of ``sign'' function. However, the convergence time is dependent on initial conditions. Then, borrowing the idea from \cite{chen20182}, two fixed-time GMs are proposed. One is designed using the property of sine function, whose convergence time is dependent on the frequency of a sine function. The other one is designed using the property of Mittag-Leffler function, whose convergence time is determined by the first positive zero of a Mittag-Leffler function. All the conclusions are extended to more general cases and finally validated by simulation examples.

The remainder of the paper is organized as follows. Section \ref{sec2} gives some basic definitions about fractional calculus and convex optimization. Finite-time GM is provided in Section \ref{sec3}. Two different types of fixed-time GM are given in Section \ref{sec4}. The paper is finally concluded in Section \ref{sec5}.

\textbf{Notations:} Throughout the paper, $f(t)*g(t)$ denotes the convolution of function $f(t)$ and $g(t)$. $\left\langle x,y\right\rangle$ implies the inner product of vector $x$ and $y$. ${\mathscr L}\{\cdot\}$ denotes the Laplace transform. $\nabla f(x)$ implies the gradient of $f(x)$. $\left\|\cdot\right\|$ indicates the Euclid norm.
\section{Preliminaries}\label{sec2}
\begin{definition}\cite{boyd2004convex}
For a function $f(x)$ whose gradient exists, if there exists a scalar $L>0$ such that
\begin{eqnarray}
\left\langle {\nabla f\left( {{\theta _1}} \right) - \nabla f\left( {{\theta _2}} \right),{\theta _1} - {\theta _2}} \right\rangle  \le L{\left\| {{\theta _1} - {\theta _2}} \right\|^2},
\end{eqnarray}
for any $\theta_1$ and $\theta_2$ belonging to the definition domain of $f(x)$, then $f(x)$ is said to have an $L$-continuous gradient.
\end{definition}
\begin{definition}\cite{boyd2004convex}
For a convex function $f(x)$ whose gradient exists, if there exists a scalar $\mu>0$ such that
\begin{eqnarray}
\left\langle {\nabla f\left( {{\theta _1}} \right) - \nabla f\left( {{\theta _2}} \right),{\theta _1} - {\theta _2}} \right\rangle  \ge \mu {\left\| {{\theta _1} - {\theta _2}} \right\|^2},
\end{eqnarray}
for any $\theta_1$ and $\theta_2$ belonging to the definition domain of $f(x)$, then $f(x)$ is said to be $\mu$-strong convex.
\end{definition}

For any constant $n-1<\alpha<n,~n\in \mathbb{N}_{+}$, the Caputo's derivative \cite{podlubny1998fractional} with order $\alpha$ for a smooth function $f(t)$ is given by
\begin{eqnarray}\label{Caputo}
{\mathscr D}^\alpha f\left( t \right) = \frac{1}{{\Gamma \left( {n - \alpha } \right)}}\int_{0}^t {\frac{{{f^{\left( n \right)}}\left( \tau  \right)}}{{{{\left( {t - \tau } \right)}^{\alpha  - n + 1}}}}{\rm{d}}\tau },
\end{eqnarray}
whose Laplace transform is
\[{\mathscr L}\left\{ {\mathscr D}^\alpha f\left( t \right) \right\} = {s^\alpha }F\left( s \right) - \sum\limits_{i = 0}^{n - 1} {{s^{\alpha  - i - 1}}{f^{\left( i \right)}}\left( 0 \right)}, \]
where $F(s)$ is the Laplace transform of $f(t)$.

Mittag-Leffler function, an extension of traditional exponential function, plays an important role in fractional calculus, which is defined as
\[{E_{\alpha ,\beta }}\left( z \right) = \sum\limits_{k = 0}^\infty  {\frac{{{z^k}}}{{\Gamma \left( {\alpha k + \beta } \right)}}},~\alpha>0,~\beta>0. \]

The following Laplace pair always holds for Mittag-Leffler function
\begin{eqnarray}\label{eqn1}
{\mathscr L}\left\{ {{t^{\beta  - 1}}{E_{\alpha ,\beta }}\left( { - \rho {t^\alpha }} \right)} \right\} = \frac{{{s^{\alpha  - \beta }}}}{{{s^\alpha } + \rho }},
\end{eqnarray}
where $\alpha>0$, $\beta>0$, and $\rho>0$.
\begin{lemma}\label{lemma1}
Mittag-Leffler function $E_{\alpha,1}(-\rho t^\alpha)$ must contain a positive zero for any $1<\alpha<2$ and $\rho>0$.
\begin{proof}
We will prove the lemma by contradiction. Consider the input function as $\frac{t^{\delta-1}}{\Gamma(\delta)},~1<\delta<\alpha$ whose Laplace transform is $\frac{1}{s^\delta}$. Then define an output as
\[y\left( t \right) = {E_{\alpha ,1 }}\left( { - \rho {t^\alpha }} \right)*\frac{{{t^{\delta  - 1}}}}{{\Gamma \left( \delta  \right)}} = {L^{ - 1}}\left\{ {\frac{{{s^{\alpha  - 1}}}}{{{s^\alpha } + \rho }}\frac{1}{{{s^\delta }}}} \right\}.\]
One can then conclude that $\mathop {\lim }\limits_{t \to \infty } y\left( t \right) = 0$ using the final value theorem of the Laplace transform.

On the other hand, assume that $E_{\alpha,1}(-\rho t^\alpha)>0$ holds for $t>0$, and one can obtain that
\[\mathop {\lim }\limits_{t \to \infty } y\left( t \right) = \int_0^\infty  {{E_{\alpha ,1 }}\left( { - \rho {{\left( {t - \tau } \right)}^\alpha }} \right)\frac{{{\tau ^{\delta  - 1}}}}{{\Gamma \left( \delta  \right)}}} {\rm d}\tau \]
must be larger than zero, which contradicts to the conclusion that $\mathop {\lim }\limits_{t \to \infty } y\left( t \right) = 0$. By contradictions, it is known that $E_{\alpha,1}(-\rho t^\alpha)$ must contain a zero for $t>0$. This completes the proof.
\end{proof}
\end{lemma}

Similarly, one has the following lemma.
\begin{lemma}\label{lemma2}
Mittag-Leffler function $E_{\alpha,\alpha}(-\rho t^\alpha)$ must contain a positive zero for any $1<\alpha<2$ and $\rho>0$.
\end{lemma}

\begin{lemma}\label{lemma3}
For any $1<\alpha<2$ and $\rho>0$, the first positive zero of $E_{\alpha,1}(-\rho t^\alpha)$ is smaller than that of $E_{\alpha,\alpha}(-\rho t^\alpha)$.
\end{lemma}
\begin{proof}
Assume the first positive zero of $E_{\alpha,1}(-\rho t^\alpha)$ is $t_0$. Since the first positive zero of $E_{\alpha,1}(-\rho t^\alpha)$ is equivalent to that of $t^{\alpha-1}E_{\alpha,1}(-\rho t^\alpha)$ whose Laplace transform is $\frac{{{1 }}}{{{s^\alpha } + \rho }}=\frac{{{s^{\alpha  - 1}}}}{{{s^\alpha } + \rho }}\frac{1}{{{s^{\alpha  - 1}}}}$. One then has that
\[{t^{\alpha  - 1}}{E_{\alpha ,\alpha }}\left( { - \rho {t^\alpha }} \right) = {E_{\alpha ,1}}\left( { - \rho {t^\alpha }} \right)*\frac{{{t^{\alpha  - 2}}}}{{\Gamma \left( {\alpha  - 1} \right)}}\]
is larger than zero for $t\le t_0$ when ${E_{\alpha ,1}}\left( { - \rho {t^\alpha }} \right)>0$. Combing with Lemma \ref{lemma2}, one has that the first positive zero of $E_{\alpha,\alpha}(-\rho t^\alpha)$ is larger than $t_0$. This competes the proof.
\end{proof}

\begin{figure}
\centering
\includegraphics[width=0.45\textwidth]{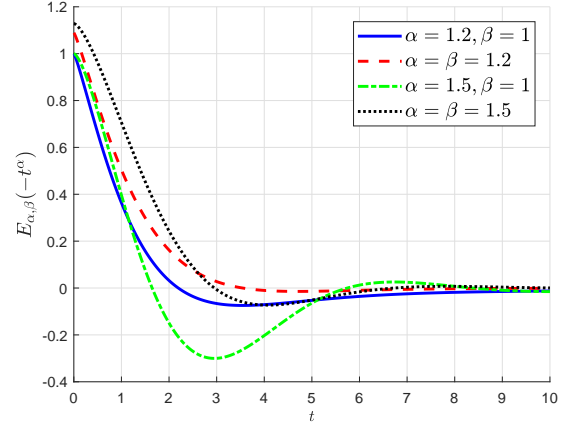}
\caption{Plots of Mittag-Leffler function with different parameter settings}\label{f1-sec2}
\end{figure}

Function plots of Mittag-Leffler function with different parameter settings are shown in Fig. \ref{f1-sec2}, and following observations can be found
\begin{itemize}
  \item Both $E_{\alpha,1}(-\rho t^\alpha)$ and $E_{\alpha,\alpha}(-\rho t^\alpha)$ with $1<\alpha<2$ have a positive zero. Moreover, a larger $\alpha$ indicates a smaller positive zero;
  \item For the same $1<\alpha<2$, the first positive zero of $E_{\alpha,1}(-\rho t^\alpha)$ is smaller than the one of $E_{\alpha,\alpha}(-\rho t^\alpha)$.
\end{itemize}
\section{Design of finite-time GM}\label{sec3}
Consider an unconstrained convex optimization problem, where the target function $f(x),~x\in {\mathbb R}^n$ is convex and has a global minimum point $x^*$. The conventional GM
\begin{equation}\label{}
\dot x =  - \rho {{\nabla f\left( x \right)}},
\end{equation}
will asymptotically converge to the minimum point if the function $f(x)$ has an $L$-continuous gradient. To achieve a finite-time convergence, one can borrow the idea from the sliding mode control and the following GM can be designed
\begin{equation}\label{ch4-552}
\dot x =  - \rho \frac{{\nabla f\left( x \right)}}{{{{\left\| {\nabla f\left( x \right)} \right\|}^2}}},
\end{equation}
where $\rho>0$ is the step size.
\begin{theorem}\label{thm1}
If the convex function $f(x)$ has an $L$-continuous gradient, then algorithm (\ref{ch4-552}) can reach the minimum point $x^*$ in a finite time.
\end{theorem}
\begin{proof}
Consider the Lyapunov function $V = {\left\| {x - {x^*}} \right\|^2}$, and take the first-order time derivative, yielding,
\[\begin{array}{rl}
\dot V =&\hspace{-6pt} 2{\left( {x - {x^*}} \right)^{\rm T}}\dot x\\
 = &\hspace{-6pt}  - 2\rho \frac{{{{\left( {x - {x^*}} \right)}^{\rm T}}\nabla f\left( x \right)}}{{{{\left\| {\nabla f\left( x \right)} \right\|}^2}}}\\
 \le &\hspace{-6pt}  - \frac{2\rho }{L}.
\end{array}\]

Take time integral on both sides, yielding,
\[\begin{array}{l}
V\left( t \right) - V\left( 0 \right) \le  - 2\frac{\rho }{L}t\\
 \Rightarrow V\left( t \right) = V\left( 0 \right) - 2\frac{\rho }{L}t.
\end{array}\]

Since $V(t)\ge 0$, one has that the minimum point can be reached in a finite time. Moreover, the convergence time is shorter than $\frac{L}{2\rho }{\left\| {{x(0)} - {x^*}} \right\|^2}$. This completes the proof.
\end{proof}

Furthermore, consider the following more general case
\begin{eqnarray}\label{ch4-111}
\dot x =  - \rho \frac{{\nabla f\left( x \right)}}{{{{\left\| {\nabla f\left( x \right)} \right\|}^\alpha }}},~0 < \alpha  < 2.
\end{eqnarray}

\begin{theorem}\label{thm2}
If the convex function $f(x)$ has an $L$-continuous gradient and is $\mu$-strong convex, then algorithm (\ref{ch4-111}) can reach the minimum point $x^*$ in a finite time.
\end{theorem}

\begin{proof}
Consider the Lyapunov function $V = {\left\| {x - {x^*}} \right\|^2}$, and take the first-order time derivative, yielding,
\begin{eqnarray}\label{11223}
\begin{array}{rl}
\dot V =&\hspace{-6pt} 2{\left( {x - {x^*}} \right)^{\rm T}}\dot x\\
 =&\hspace{-6pt}  - 2\rho \frac{{{{\left( {x - {x^*}} \right)}^{\rm T}}\nabla f\left( x \right)}}{{{{\left\| {\nabla f\left( x \right)} \right\|}^\alpha }}}\\
 \le &\hspace{-6pt} - \frac{2\rho }{L}{\left\| {\nabla f\left( x \right)} \right\|^{2 - \alpha }}\\
 \le &\hspace{-6pt} - \frac{{2\rho {\mu ^{2-\alpha}}}}{L}{\left\| {x - {x^*}} \right\|^{2 - \alpha }}\\
 = &\hspace{-6pt} - \frac{{2\rho {\mu ^{2-\alpha}}}}{L}{V^{\frac{{2 - \alpha }}{2}}},
\end{array}
\end{eqnarray}
where the $L$-continuous gradient and strong convex properties are used.

Take time integral on both sides of (\ref{11223}), yielding,
\[\begin{array}{l}
\frac{2}{\alpha }{V^{\frac{\alpha }{2}}}\left( t \right) - \frac{2}{\alpha }{V^{\frac{\alpha }{2}}}\left( 0 \right) \le  - \frac{{2\rho {\mu ^{2 - \alpha }}}}{L}t\\
 \Rightarrow V\left( t \right) \le {\left[ {{V^{\frac{\alpha }{2}}}\left( 0 \right) - \frac{{\rho {\mu ^{2 - \alpha }}}}{{\alpha L}}t} \right]^{{2 \mathord{\left/
 {\vphantom {2 \alpha }} \right.
 \kern-\nulldelimiterspace} \alpha }}}.
\end{array}\]

Since $V(t)\ge 0$, one has that the minimum point can be reached in a finite time. Moreover, the convergence time is shorter than $\frac{{L}}{{\rho {\mu ^{2-\alpha}}\alpha }}{\left\| {{x(0)} - {x^*}} \right\|^\alpha }$. This completes the proof.
\end{proof}
\begin{remark}\label{re1}
Here, some comments on Theorem \ref{thm1} and Theorem \ref{thm2} are given.
\begin{itemize}
  \item When $\alpha=2$, Theorem \ref{thm2} will reduce to Theorem \ref{thm1}, and the condition ``$\mu$-strong convex'' is not necessary any more.
  \item For both algorithms (\ref{ch4-552}) and (\ref{ch4-111}), the convergence time is sensitive to initial conditions since the upper bounds directly contain the initial value $x(0)$.
  \item If $\alpha$ tends to zero, then the upper bound of convergence time using algorithm (\ref{ch4-111}) will tend to infinity, which is reduced to the traditional GM and results in an asymptotic convergence.
  \item Both algorithms (\ref{ch4-552}) and (\ref{ch4-111}) will be singular when $\left\|\nabla f(x)\right\|=0$. Therefore, one can add a small positive scalar to avoid the singularity for practical usage. Then algorithm (\ref{ch4-111}) can be modified as
      \begin{eqnarray}\label{ch4-114}
       \dot x =  - \rho \frac{{\nabla f\left( x \right)}}{{{{\left(\left\| {\nabla f\left( x \right)} \right\|+\delta\right)}^\alpha }}},~0 < \alpha  < 2,~\delta>0.
       \end{eqnarray}
\end{itemize}
\end{remark}
\begin{example}\label{ex1}
Consider the following convex function which has an $L$-continuous gradient ($L=4.30$) and is $\mu$-strong convex ($\mu=0.70$)
\[f\left( x \right) = \left[ {{x_1},{x_2}} \right]\left[ {\begin{array}{*{20}{c}}
1&1\\
1&4
\end{array}} \right]\left[ {\begin{array}{*{20}{c}}
{{x_1}}\\
{{x_2}}
\end{array}} \right].\]
Take step size $\rho=10$, initial value $x(0)=[10,-10]^{\rm T}$, and $\delta=0.01$ when simulating. For different order $\alpha$, the results are shown in Fig. \ref{f1}. For different order $\alpha$, finite-time convergence can be realized. Moreover, a smaller $\alpha$ results in a shorter convergence time.

\begin{figure}
\centering
\includegraphics[width=0.45\textwidth]{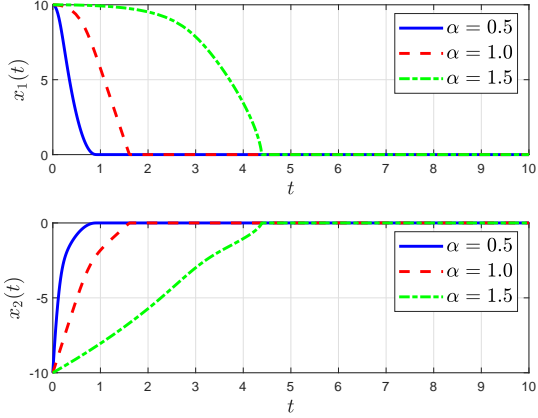}
\caption{Convergence results with different $\alpha$ in Example \ref{ex1}}\label{f1}
\end{figure}

For different initial values with fixed order $\alpha=1.0$, results are shown in Fig. \ref{f2}. It is found that larger ${\left\| {{x(0)} - {x^*}} \right\| }$ leads to longer convergent time. Moreover, the convergence time of different state variables is almost the same from both Fig. \ref{f1} and Fig. \ref{f2}, which indicates that the convergence speed for different states is robust to the ``condition number''.
\begin{figure}
\centering
\includegraphics[width=0.45\textwidth]{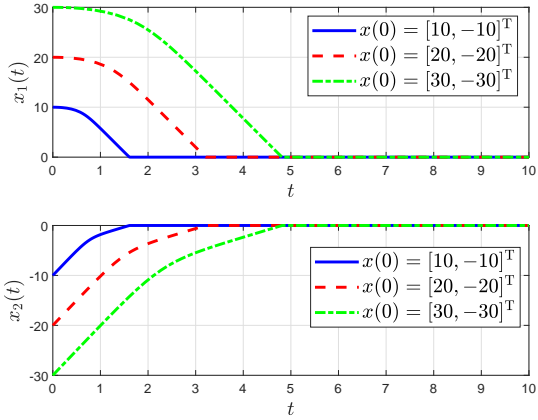}
\caption{Convergence results with different $x(0)$ in Example \ref{ex1}}\label{f2}
\end{figure}

\end{example}
\section{Design of fixed-time GM}\label{sec4}
The convergence time of algorithm (\ref{ch4-111}) is dependent on initial conditions, and can be arbitrarily large as the initial condition grows, which is undesirable. To make the convergence time robust to initial conditions, two different types of fixed-time GMs are designed in this section.
\subsection{Second-order design}
In this subsection, the fixed-time convergence is realized by using the property of sine function. The following GM can be designed
\begin{eqnarray}\label{ch4-113}
\left\{ \begin{array}{l}
\dot x =  - \theta \frac{{\nabla f\left( x \right)}}{{{{\left\| {\nabla f\left( x \right)} \right\|}^2}}},\\
\dot \theta  =  - \lambda \theta  + \rho {\left\| {\nabla f\left( x \right)} \right\|^2},
\end{array} \right.
\end{eqnarray}
where $\lambda>0,~\rho>0$.

\begin{theorem}\label{thm3}
For any convex function $f(x)$ which has an $L$-continuous gradient and is $\mu$-strong convex, GM (\ref{ch4-113}) can reach the minimum point in a finite time if the condition ${\lambda ^2} < \frac{{8\rho {\mu ^2}}}{L}$ is satisfied. Moreover, the convergence time is no longer than $\frac{\pi}{\sqrt {\frac{{4\rho {\mu ^2}}}{L} - \frac{{{\lambda ^2}}}{4}}}$.
\end{theorem}
\begin{proof}
Consider the Lyapunov function as $V = {\left\| {x - {x^*}} \right\|^2}$. Take the first-order time derivative and use the condition of L-continuous gradient, yielding,
\[\dot V =  - \theta \frac{{2{{\left( {x - {x^*}} \right)}^{\rm T}}\nabla f\left( x \right)}}{{\left\| {\nabla f\left( x \right)} \right\|^2}} \le  - \frac{{2\theta }}{L}.\]
According to the second equality in (\ref{ch4-113}), one has that $\theta(t)\ge 0$ since $\theta(0)=0$ and $\left\| {\nabla f\left( x \right)} \right\|\ge 0$. Then, one has that $\dot V\le 0$. Moreover, function $f(x)$ is $\mu$-strong convex, which indicates
\[\dot \theta  =  - \lambda \theta  + \rho {\left\| {\nabla f\left( x \right)} \right\|^2} \ge  - \lambda \theta  + \rho {\mu ^2}V.\]

Therefore, the following equalities hold
\[\left\{ \begin{array}{l}
\dot V =  - \frac{2\theta }{L} - M\left( t \right),\\
\dot \theta  =  - \lambda \theta  + \rho {\mu ^2}V + N\left( t \right),
\end{array} \right.\]
where $M(t)\ge0$ and $N(t)\ge 0$ hold for any $t\ge 0$. Take Laplace transform on both sides, yielding,
\begin{eqnarray}\label{ch4-115}
\left\{ \begin{array}{l}
s V\left( s \right) - {V(0)} =  - \frac{2}{L}\Theta \left( s \right) - M\left( s \right),\\
s\Theta \left( s \right) =  - \lambda \Theta \left( s \right) + \rho {\mu ^2} V\left( s \right) +  N\left( s \right),
\end{array} \right.
\end{eqnarray}
where $V(s)$, $\Theta (s)$, $M(s)$, and $N(s)$ are the corresponding Laplace transform of $V(t)$, $\theta (t)$, $M(t)$, and $N(t)$. $V(s)$ can then be derived by solving (\ref{ch4-115})
\begin{eqnarray}\label{ch4-120}
\begin{array}{rl}
V\left( s \right) =&\hspace{-6pt} \frac{{s + \lambda }}{{{s^2} + \lambda s + \frac{{2\rho {\mu ^2}}}{L}}}\left( {V(0) - M\left( s \right)} \right)\\
 &\hspace{-6pt}- \frac{{\frac{2}{L}}}{{{s^2} + \lambda s + \frac{{2\rho {\mu ^2}}}{L}}}N\left( s \right).
\end{array}
\end{eqnarray}

Since $\dot V \le 0$ and $V(t)\ge 0$, we only need to prove that $V(t)$ will reach zero in a finite time to achieve a finite-time convergence. To simplify the expression, define
\[\omega : = \sqrt {\frac{{2\rho {\mu ^2}}}{L} - \frac{{{\lambda ^2}}}{4}} ,\]
which is a positive real number according to the condition ${\lambda ^2} < \frac{{8\rho {\mu ^2}}}{L}$.

Perform inverse Laplace transform on both sides of (\ref{ch4-120}), resulting in
\[\begin{array}{rl}
V\left( t \right) =&\hspace{-6pt} {{\rm e}^{ - \frac{\lambda }{2}t}}\left( {\cos \left( {\omega t} \right) + \frac{\lambda }{{\omega }}\sin \left( {\omega t} \right)} \right)*\left( {V(0) - M\left( t \right)} \right)\\
 &\hspace{-6pt}- \frac{{2{{\rm e}^{ - \frac{\lambda }{2}t}}}}{{L\omega }}\sin \left( {\omega t} \right)*N\left( t \right).
\end{array}\]

The first positive zero $t_0$ of function $f(t)={\cos \left( {\omega t} \right) + \frac{\lambda }{{2\omega }}\sin \left( {\omega t} \right)}$ must be smaller than $\pi/\omega$. Thus, for any $0<t\le t_0\le \pi/\omega$, one has $f(t)\ge 0$ and $\sin(\omega t)>0$. Combining with $M(t)\ge 0$ and $N(t)\ge 0$, following inequalities hold
\[ - {{\rm e}^{ - \frac{\lambda }{2}t}}\left( {\cos \left( {\omega t} \right) + \frac{\lambda }{{\omega }}\sin \left( {\omega t} \right)} \right)*M\left( t \right)<0,\]
and
\[ - \frac{2}{{L\omega }}{{\rm e}^{ - \frac{\lambda }{2}t}}\sin \left( {\omega t} \right)*N\left( t \right)<0.\]

Then the following inequality can be derived
\[V\left( t \right) < {{\rm e}^{ - \frac{\lambda }{2}t}}\left( {\cos \left( {\omega t} \right) + \frac{\lambda }{{\omega }}\sin \left( {\omega t} \right)} \right){V(0)}.\]

Since function $f(t)={\cos \left( {\omega t} \right) + \frac{\lambda }{{2\omega }}\sin \left( {\omega t} \right)}$ must have a zero in half cycle, then $V(t)$ must reach zero within $\pi/\omega$. Combing with that $\dot V (t)\le 0$ and $V(t)\ge 0$, it is known that $V(t)$ will reach and stay on zero in a finite time. Moreover, the convergence time is shorter than $\pi/\omega$, which has no relation to initial conditions and thus indicates a fixed-time convergence. This completes the proof.
\end{proof}

Furthermore, consider the following more general case
\begin{eqnarray}\label{ch4-118}
\left\{ \begin{array}{l}
\dot x =  - \theta \frac{{\nabla f\left( x \right)}}{{{{\left\| {\nabla f\left( x \right)} \right\|}^\alpha }}},\\
\dot \theta  =  - \lambda \theta  + \rho {\left\| {\nabla f\left( x \right)} \right\|^\alpha },
\end{array} \right.
\end{eqnarray}
where $\lambda>0,~\rho>0,~0<\alpha<2$.

\begin{theorem}\label{thm4}
For any convex function $f(x)$ which has an $L$-continuous gradient and is $\mu$-strong convex, GM (\ref{ch4-113}) can reach the minimum point in a finite time if the condition ${\lambda ^2} < \frac{{8\rho {\mu ^2}}}{\alpha L}$ is satisfied. Moreover, the convergence time is no longer than $\frac{\pi}{\sqrt {\frac{{4\rho {\mu ^2}}}{\alpha L} - \frac{{{\lambda ^2}}}{4}}}$.
\end{theorem}
\begin{proof}
Take the Lyapunov function as $V = {\left\| {x - {x^*}} \right\|^2}$, and take the first-order time derivative, yielding,
\[\begin{array}{rl}
\dot V = &\hspace{-6pt} - \theta \frac{{2{{\left( {x - {x^*}} \right)}^{\rm T}}\nabla f\left( x \right)}}{{{{\left\| {\nabla f\left( x \right)} \right\|}^\alpha }}}\\
 \le&\hspace{-6pt}  - \frac{{2\theta }}{L}{\left\| {\nabla f\left( x \right)} \right\|^{2 - \alpha }}\\
 \le&\hspace{-6pt}  - \frac{{2{\mu ^{2 - \alpha }}\theta }}{L}{\left\| {x - {x^*}} \right\|^{2 - \alpha }}\\
 =&\hspace{-6pt}  - \frac{{2{\mu ^{2 - \alpha }}\theta }}{L}{V^{\frac{{2 - \alpha }}{2}}}.
\end{array}\]

Defining $\hat V := {V^{\frac{\alpha }{2}}}$, one has that
\[\dot {\hat V} \le  - \frac{{4 {\mu ^{2 - \alpha }}}}{\alpha L}\theta ,\]
and
\[\dot \theta  =  - \lambda \theta  + \rho {\left\| {\nabla f\left( x \right)} \right\|^\alpha } \ge  - \lambda \theta  + \rho {\mu ^\alpha }\hat V.\]

Similar to the proof of Theorem \ref{thm3}, the proof can be completed.
\end{proof}
\begin{remark}
Some comments on Theorem \ref{thm3} and Theorem \ref{thm4} are given as follows.
\begin{itemize}
  \item Theorem \ref{thm4} will reduce to Theorem \ref{thm3} when $\alpha=2$. Moreover, algorithm (\ref{ch4-118}) with $\alpha=0$ will finally indicate an asymptotic convergence.
  \item Since the upper bound for the convergence time is determined by $\frac{\pi}{\sqrt {\frac{{4\rho {\mu ^2}}}{\alpha L} - \frac{{{\lambda ^2}}}{4}}}$, one can set $\alpha=1.0$ and tune $\lambda$ to achieve a desirable convergence time for practical usage.
  \item The parameter $\lambda$ in algorithm (\ref{ch4-118}) is used to attenuate the value of $\theta$ after reaching the minimum point. It is quite useful when realizing algorithm (\ref{ch4-118}) in its discretization form. When $\lambda=0$, the attenuating item ${\rm e}^{-\frac{\lambda}{2}t}$ will disappear during the proof process of Theorem \ref{thm3} and the conclusion for fixed-time convergence still holds.
  \item To avoid singularity, an additional positive scalar can be introduced for practical usage and algorithm (\ref{ch4-118}) can be modified as
      \begin{eqnarray}\label{ch4-1181}
\left\{ \begin{array}{l}
\dot x =  - \theta \frac{{\nabla f\left( x \right)}}{{{{\left(\left\| {\nabla f\left( x \right)} \right\|+\delta\right)}^\alpha }}},\\
\dot \theta  =  - \lambda \theta  + \rho {\left\| {\nabla f\left( x \right)} \right\|^\alpha },
\end{array} \right.
\end{eqnarray}
where $\delta>0$ is a small scalar.
\end{itemize}
\end{remark}
\begin{example}\label{ex2}
Consider the same convex function in Example \ref{ex1}. Set step size $\rho=10$, $\lambda=1$, and $\delta=0.01$ when simulating. For different initial values with the same order $\alpha=1$, the results are shown in Fig. \ref{f3}. The upper bound is estimated as $\frac{\pi}{\sqrt {\frac{{4\rho {\mu ^2}}}{\alpha L} - \frac{{{\lambda ^2}}}{4}}}$, i.e., $1.51$ {\rm{(sec)}}. It is observed that algorithm (\ref{ch4-1181}) reaches the minimum point at almost the same time (about $0.45$ {\rm{sec}}, smaller than the estimated upper bound) from different initial values, which demonstrates the results in Theorem \ref{thm4}.

\begin{figure}
\centering
\includegraphics[width=0.45\textwidth]{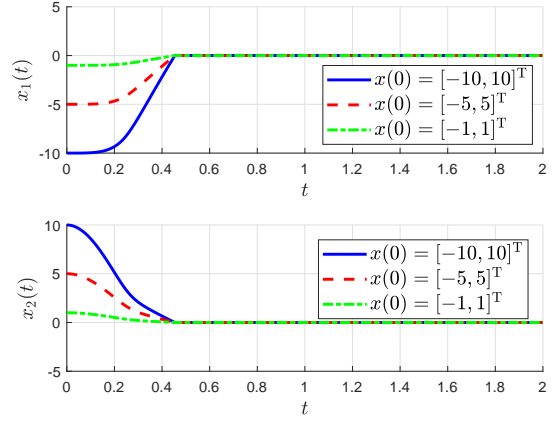}
\caption{Convergence results with different $x(0)$ in Example \ref{ex2}}\label{f3}
\end{figure}

For different $\alpha$ with the same initial value $x(0)=[-5,5]^{\rm T}$, results are shown in Fig. \ref{f4}. It is found that fixed-time convergence can be reached in all cases. Moreover, order $\alpha$ only has a small influence on the convergence time, thus one can choose $\alpha=1.0$ for practical usage.

\begin{figure}
\centering
\includegraphics[width=0.45\textwidth]{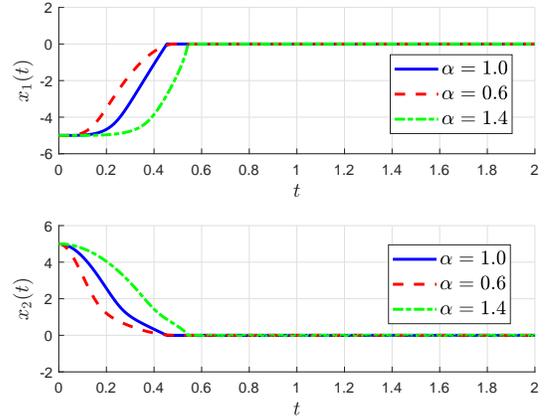}
\caption{Convergence results with different $\alpha$ in Example \ref{ex2}}\label{f4}
\end{figure}

\end{example}
\begin{example}\label{ex4}
In this example, Zakharov function is considered, which can be formulated as
\[f\left( x \right) = \sum\limits_{i = 1}^n {x_i^2}  + {\left( {\sum\limits_{i = 1}^n {0.5i{x_i}} } \right)^2} + {\left( {\sum\limits_{i = 1}^n {0.5i{x_i}} } \right)^4},\]
which is not strong convex nor has an $L$-gradient continuous globally. When simulating, set $n=2$, $\alpha=1.0$, $\rho=10$, and $\delta=0.01$. Results with different initial conditions are shown in Fig. \ref{f10}. It is found that finite-time convergence can be obtained while fixed-time convergence cannot be guaranteed any more. Generally, a larger initial condition will lead to a shorter convergence time. As $x\to 0$, the quadratic terms in Zakharov function become dominant, and Zakharov function will then has an $L$-continuous gradient. The finite-time convergence can be derived according to Theorem \ref{thm2}.

Results with different order $\alpha$ are shown in Fig. \ref{f11}. It is found that finite-time convergence can still be derived while the convergence time is more sensitive to order varying compared with the results in Example \ref{ex2}. Moreover, a larger $\alpha$ leads to a shorter convergence time in this example.

\begin{figure}
\centering
\includegraphics[width=0.45\textwidth]{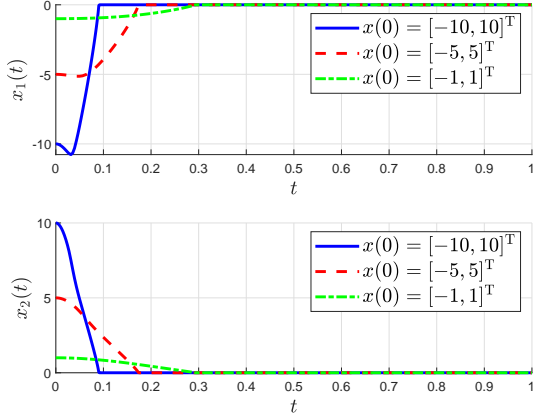}
\caption{Convergence results with different initial conditions in Example \ref{ex4}}\label{f11}
\end{figure}
\begin{figure}
\centering
\includegraphics[width=0.45\textwidth]{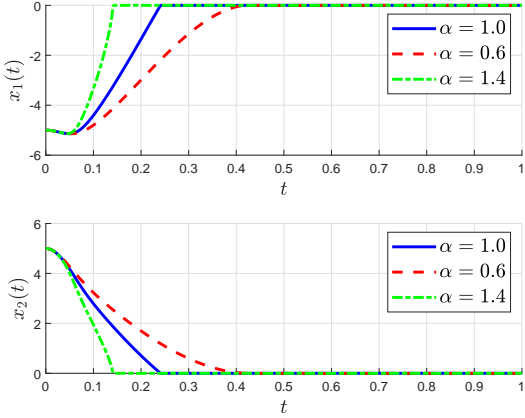}
\caption{Convergence results with different order $\alpha$ in Example \ref{ex4}}\label{f10}
\end{figure}
\end{example}
\subsection{Fractional-order design}
Fractional calculus is a natural extension of integer order calculus, which plays an important role in all kinds of fields. In this subsection, fractional design of finite-time GM is considered by replacing the update law for $\theta$ in algorithm (\ref{ch4-113}) with a fractional one, which can be formulated as
\begin{eqnarray}\label{ch4-222}
\left\{ \begin{array}{l}
\dot x =  - \theta \frac{{\nabla f\left( x \right)}}{{{{\left\| {\nabla f\left( x \right)} \right\|}^2}}}\\
{\mathscr D}^\beta \theta  = \rho {\left\| {\nabla f\left( x \right)} \right\|^2}
\end{array} \right.,
\end{eqnarray}
where $\theta(0)=0$ and $0<\beta<1$.

\begin{theorem}\label{thm5}
For any convex function $f(x)$ which has an $L$-continuous gradient and is $\mu$-strong convex, algorithm (\ref{ch4-222}) can reach the minimum point in a finite time. Moreover, the convergence time is shorter than the first positive zero of Mittag-Leffler function ${E_{\beta  + 1,1}}\left( { - \frac{2\rho {\mu ^2}}{L}{t^{ {1 + \beta } }}} \right)$.
\end{theorem}
\begin{proof}
Similar to the proof of Theorem \ref{thm3}. Consider the Lyapunov function as $V = {\left\| {x - {x^*}} \right\|^2}$, and use the condition of $L$-continuous gradient, yielding,
\[\dot V =  - \theta \frac{{2{{\left( {x - {x^*}} \right)}^{\rm T}}\nabla f\left( x \right)}}{{\left\| {\nabla f\left( x \right)} \right\|^2}} \le  - \frac{{2\theta }}{L}.\]
Combining with that $\theta(t)\ge 0$, one has that $\dot V\le 0$. Moreover, function $f(x)$ is $\mu$-strong convex, which indicates
\[{{\mathscr D}^\beta }\theta  \ge \rho {\mu ^2}V.\]

Then, one can arrive at following equalities
\begin{eqnarray}\label{eq1123}
\left\{ \begin{array}{l}
\dot V =  - \frac{{2\theta }}{L} - M\left( t \right),\\
{{\mathscr D}^\beta }\theta  = \rho {\mu ^2}V + N\left( t \right),
\end{array} \right.
\end{eqnarray}
where $M(t)\ge 0$ and $N(t)\ge 0$. Perform Laplace transform on both sides of (\ref{eq1123}), resulting in
\begin{eqnarray}\label{eq123}
\left\{ \begin{array}{l}
sV\left( s \right) - {V(0)} =  - \frac{2}{L}\Theta \left( s \right) - M\left( s \right),\\
{s^\beta }\Theta \left( s \right) = \rho {\mu ^2}V\left( s \right) + N\left( s \right),
\end{array} \right.
\end{eqnarray}
where $\theta(0)=0$ is used. $V(s)$ can then be derived by solving (\ref{eq123})
\begin{eqnarray}\label{ch4-555}
V\left( s \right) = \frac{{{s^\beta }}}{{{s^{\beta  + 1}} + \frac{2\rho {\mu ^2}}{L}}}{V(0)} - \frac{{{s^\beta }M\left( s \right) + \frac{2}{L}N\left( s \right)}}{{{s^{\beta  + 1}} + \frac{2\rho {\mu ^2}}{L}}}.
\end{eqnarray}

Perform inverse Laplace transform on both sides of (\ref{ch4-555}), yielding,
\[\begin{array}{rl}
V\left( t \right) =&\hspace{-6pt} {E_{\beta  + 1,1}}\left( { - \frac{2\rho {\mu ^2}}{L}{t^{ {1 + \beta } }}} \right)*\left( {V(0) - M\left( t \right)} \right)\\
 &\hspace{-6pt}- {t^\beta }{E_{\beta  + 1,\beta  + 1}}\left( { - \frac{2\rho {\mu ^2}}{L}{t^{ {1 + \beta } }}} \right)*N\left( t \right).
\end{array}\]

Suppose $t_0$ and $t_1$ are the corresponding first positive zero of ${E_{\beta  + 1,1}}\left( { - \frac{2\rho {\mu ^2}}{L}{t^{ {1 + \beta } }}} \right)$ and ${t^\beta E_{\beta  + 1,\beta+1}}\left( { - \frac{2\rho {\mu ^2}}{L}{t^{ {1 + \beta } }}} \right)$. According to Lemma \ref{lemma3}, one has that $t_0<t_1$. Then for any $0<t<t_0$, the following two inequalities hold,
\[{E_{\beta  + 1,1}}\left( { - \frac{2\rho {\mu ^2}}{L}{t^{ {1 + \beta } }}} \right)>0,\]
and
\[{t^\beta E_{\beta  + 1,\beta+1}}\left( { - \frac{2\rho {\mu ^2}}{L}{t^{ {1 + \beta } }}} \right)>0.\]

Finally, we arrive at the following inequality
\[V\left( t \right) \le {V(0)}{E_{\beta  + 1,1}}\left( { \frac{2\rho {\mu ^2}}{L}{t^{ {1 + \beta } }}} \right).\]

On one hand, function ${E_{\beta  + 1,1}}\left( { \frac{2\rho {\mu ^2}}{L}{t^{ {1 + \beta } }}} \right)$ must have a positive zero for $0<\beta<1$, thus $V(t)$ must reach its zero in a finite time. On the other hand, $V(t)$ will maintain on zero once $V(t)$ reaches zero since $\dot V(t)\le 0$. Moreover, the first positive zero of ${E_{\beta  + 1,1}}\left( { - \frac{2\rho {\mu ^2}}{L}{t^{ {1 + \beta } }}} \right)$ is the upper bound of the convergence time, which is irrelevant to initial conditions and thus indicates a fixed-time convergence. This completes the proof.
\end{proof}

Furthermore, consider the following general case
\begin{eqnarray}\label{ch4-2222}
\left\{ \begin{array}{l}
\dot x =  - \theta \frac{{\nabla f\left( x \right)}}{{{{\left\| {\nabla f\left( x \right)} \right\|}^\alpha }}}\\
{\mathscr D}^\beta \theta  = \rho {\left\| {\nabla f\left( x \right)} \right\|^\alpha }
\end{array} \right.,~0<\alpha<2,~0<\beta<1.
\end{eqnarray}
\begin{theorem}\label{thm6}
For any convex function $f(x)$ which has an $L$-continuous gradient and is $\mu$-strong convex, algorithm (\ref{ch4-2222}) can reach the minimum point in a finite time. Moreover, the convergence time is shorter than the first positive zero of Mittag-Leffler function ${E_{\beta  + 1,1}}\left( { - \frac{4\rho\mu^2}{\alpha L}{t^{ {1 + \beta } }}} \right)$.
\end{theorem}
\begin{proof}
Similar to the proof of Theorem \ref{thm4}, consider the Lyapunov function $V = {\left\| {x - {x^*}} \right\|^2}$ and one has
\begin{eqnarray}\label{eq1124}
\left\{ \begin{array}{l}
\dot {\hat V} =  - \frac{{4 {\mu ^{2 - \alpha }}}}{\alpha L}\theta - M\left( t \right),\\
{{\mathscr D}^\beta }\theta  = \rho {\mu ^\alpha}\hat V + N\left( t \right),
\end{array} \right.
\end{eqnarray}
where $\hat V = {V^{\frac{\alpha }{2}}}$, $M(t)>0$, and $N(t)>0$.

Then similar to the proof of Theorem \ref{thm5}, one has that algorithm (\ref{ch4-2222}) can reach the minimum point in a fixed time, which is shorter than the first positive zero of function ${E_{\beta  + 1,1}}\left( { - \frac{4\rho\mu^2}{\alpha L}{t^{ {1 + \beta } }}} \right)$. This completes the proof. (When $\alpha=2$, Theorem \ref{thm6} reduces to Theorem \ref{thm5})
\end{proof}
\begin{remark}\label{remark4}
According to the results of Theorem \ref{thm5} and Theorem \ref{thm6}, the upper bound of the convergence time is determined by the first positive zero of a Mittag-Leffler function. Therefore, some numerical results about the first positive zero of the standard Mittag-Leffler function $E_{\beta,1}(-t^\beta)$ is provided in TABLE \ref{tab1} to help designing the parameters. Since the first positive zero of the standard Mittag-Leffler function is fixed for some specific $\beta$, parameter $\rho$, $\alpha$ and $\beta$ in algorithm (\ref{ch4-2222}) will influence the position of the first positive zero. It is straightforward that a larger $\rho$, a smaller $\alpha$, and a larger $\beta$ all indicate a smaller positive zero, which results in a smaller upper bound of fixed convergence time.
\begin{table}[!htbp]
\begin{center}
\caption{First positive zero of Mittag-Leffler function $E_{\beta,1}(-t^\beta)$}\label{tab1}
\begin{tabular}{cccccc}
\hline
$\beta$ & 1.7 & 1.5 & 1.3 & 1.1 &1.05 \\\hline
$t_{\beta}$ &1.57&1.65&1.89&2.88& 3.72\\\hline
\end{tabular}
\end{center}
\end{table}
\end{remark}
\begin{example}\label{ex3}
Consider the same convex function in Example \ref{ex1}. Take $\lambda=1$, $\alpha=1.0$. Fig. \ref{f5} shows the results of different initial values with $\beta=0.2$ and $\rho=10$. Fixed-time convergence can be directly observed. Interestingly, the larger initial condition comes into a shorter convergence time in this example.
\begin{figure}
\centering
\includegraphics[width=0.45\textwidth]{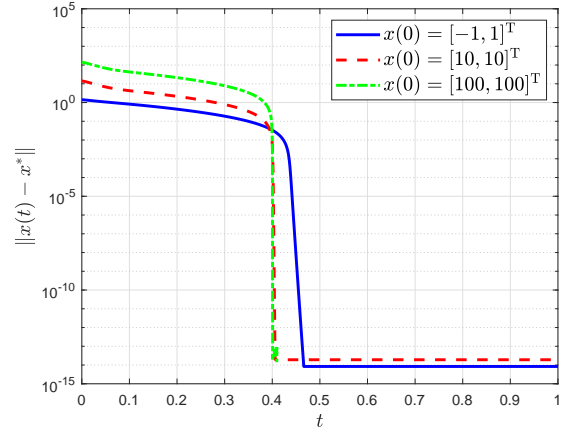}
\caption{Convergence results with different $x(0)$ in Example \ref{ex3}}\label{f5}
\end{figure}

Fig. \ref{f6} shows the results of different order $\beta$ with $x(0)=[-10,10]^{\rm T}$ and $\rho=10$. Fixed-time convergence can be directly observed for different $\beta$. Moreover, a larger $\beta$ means a shorter fixed convergence time as declared in Remark \ref{remark4}.

\begin{figure}
\centering
\includegraphics[width=0.45\textwidth]{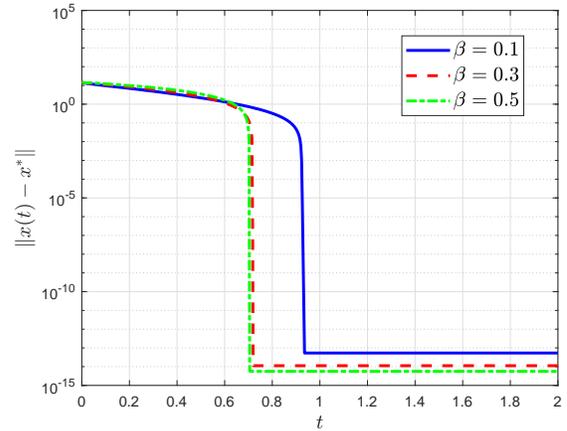}
\caption{Convergence results with different $\beta$ in Example \ref{ex3}}\label{f6}
\end{figure}

\end{example}
\section{Conclusion}\label{sec5}
In this paper, by the advantage of finite-time convergence in system control, several novel GMs have been proposed to realize finite-time and fixed-time convergence rather than asymptotic convergence. At first, finite-time convergence is derived by normalizing the gradient, but the convergence time is dependent on initial conditions. To make the convergence time robust to initial conditions, two fixed-time GMs are then provided by using the property of periodic function and Mittag-Leffler function respectively. All these results are extended to more general cases and finally demonstrated by numerical examples. There are some promising directions for future research:
\begin{itemize}
  \item designing finite-time and fixed-time GMs for a more general class of functions, such as non-strong convex functions;
  \item designing finite-time and fixed-time GMs in their discrete forms, which will be more potable for practical usage.
\end{itemize}

\bibliographystyle{IEEEtran}
\bibliography{IEEEabrv,FOEP}







\end{document}